\documentclass[12pt,reqno]{amsart}
\usepackage{verbatim}
\usepackage{bm}
\usepackage{amsfonts}
\usepackage{amsmath}
\usepackage{amssymb}
\usepackage[all]{xy}
\newtheorem{theorem}{Theorem}
\newtheorem{proposition}{Proposition}
\newtheorem{corollary}{Corollary}

\newtheorem{lemma}{Lemma}

\newfont{\bb}{msbm10 at 12pt}

\def\qed{\hfill{Q.E.D.}\smallskip}

\newcommand{\ls}{\setlength{\baselineskip}{12pt}
                 \setlength{\parskip}{3mm}}

\newcommand{\mysection}[1]{\section{#1}\setcounter{equation}{0}}

\newcommand{\bal}{\begin{align}}      \newcommand{\eal}{\end{align}}
\newcommand{\ba}{\begin{array}}      \newcommand{\ea}{\end{array}}
\newcommand{\bc}{\begin{center}}     \newcommand{\ec}{\end{center}}
\newcommand{\be}{\begin{enumerate}}  \newcommand{\ee}{\end{enumerate}}
\newcommand{\beq}{\begin{eqnarray}}  \newcommand{\eeq}{\end{eqnarray}}
\newcommand{\beQ}{\begin{eqnarray*}} \newcommand{\eeQ}{\end{eqnarray*}}
\newcommand{\bi}{\begin{itemize}}    \newcommand{\ei}{\end{itemize}}
\newcommand{\bt}{\begin{tabular}}    \newcommand{\et}{\end{tabular}}
\newcommand{\bdm}{\begin{displaymath}} \newcommand{\edm}{\end{displaymath}}

\begin{document}

\title{  Eigenvalue estimate for the Dirac-Witten operator on locally reducible Riemannian manifolds }

{\address{Yongfa Chen
 School of Mathematics and Statistics,
Central China Normal University, Wuhan 430079, P.R.China.}

\author{Yongfa Chen}}

\thanks{ }

\email{yfchen@mail.ccnu.edu.cn}

\begin{abstract}
We get optimal lower bounds for the eigenvalues of  the Dirac-Witten operator on locally reducible spacelike submanifold
 in terms of intrinsic and extrinsic expressions. The limiting-cases are also studied.
\end{abstract}

\keywords{  Dirac operator, eigenvalue,  mean curvature,
scalar curvature, reducible manifold}

\subjclass[2008]{}

\maketitle

\mysection{\textbf{Introduction}} \ls
    It is well known  that  the spectrum of
the Dirac operator on closed  spin manifolds detects
subtle information on the geometry and the topology of such
manifolds (see \cite{LM}).
A fundamental tool to get estimates for eigenvalues of  the basic Dirac operator $D$ acting on spinors
is the Schr\"{o}dinger-Lichnerowicz  formula
\beq \label{S-L1}
D^2=\nabla^*\nabla+\frac{1}{4}{R},
\eeq
where $\nabla^*$ is the formal adjoint of $\nabla$  with respect to the natural Hermitian
inner product on spinor bundle $\Sigma {M^m}$ and  $R^M$ stands for  the
scalar curvature of  the closed spin manifold $(M^m,g)$. Using (\ref{S-L1}) and a modified spin connection, Friedrich proved  any eigenvalue $\lambda$ of $D$ satisfies the following sharp inequality:
 \beq\label{Friedrich}
 \lambda^2\geqslant \frac{m}{4(m-1)} \min_M R.
\eeq
The case of equality in (\ref{Friedrich}) occurs iff  $(M^m,g)$  admits   a real {Killing spinor}, which by definition
satisfies
 the  overdetermined elliptic equation
$
\nabla_X\psi=-\frac{\lambda}{m}X\cdot \psi,
$ here $\nabla_X$ denotes the Levi-Civita connection of $g.$
The manifold  must be a locally irreducible Einstein manifold. Note complete simply-connected Riemannian spin manifolds $(M^m, g)$ carrying a
non-zero space of real Killing spinors have been completely classified by B\"ar\cite{Bar93}.

Now we shall consider another important Dirac-type operator in the field of mathematical physics. Suppose $(M^m,g)$ is the spacelike hypersurface  in an $(n+1)$-dimensional Lorentzian  $(\tilde{M}^{m,1},\tilde{g}),$
with the induced Riemannian metric $g$ and timelike unit normal vector $e_0.$
Here $\tilde{g}$  satisfies the Einstein field equations
$$ \widetilde{Ric}-\frac{1}{2}\tilde{R}\tilde{g} =T,$$
where $\widetilde{Ric},\tilde{R}$ are the Ricci curvature, scalar
curvature of $\tilde{g}$ respectively, and $T$ is the
energy-momentum tensor. We also denote the Levi-Civita connection of $\tilde{g}$ by $\tilde{\nabla}.$ Choosing an orthonormal frame
$\{e_\alpha\}$ with ${e_0}$ timelike and $\{e_i\}_{i=1}^m  $ spacelike. Then, in physics $T_{00}$ is
interpreted as
the local mass density, and $T_{0i}$ is interpreted as the local momentum density.
The well-known spinorial proof of the positive
mass theorem given by  Witten \cite{Witten} is based on  the following
Weitzenb\"{o}ck type formula
\beq\label{Wittenformula}
\tilde{D}^2=\tilde{\nabla}^*\tilde{\nabla}+\frac{1}{2}(T_{00}+T_{0i}e^0\cdot e^i\cdot),
 \eeq
 for
the Dirac-Witten operator $\tilde{D}:=e^i\cdot \tilde{\nabla}_i. $ The mass is given by the limit
at infinity of some
 boundary integral term. In view of the importance of the Dirac-Witten operator, by (\ref{Wittenformula})  Hijazi and Zhang  first established the optimal lower
bounds  for the eigenvalues of Dirac-Witten operator of compact (with or without boundary) spacelike hypersurfaces of a Lorentzian manifold,
 whose  metric satisfies the Einstein field equations and whose energy-momentum tensor satisfies
  the dominant energy condition \cite{HZ03}.
The limiting cases were also
studied.
Note,  in term of Gauss and Codazzi equation of the spacelike hypersurface $M^m\hookrightarrow \tilde{M}^{m,1}$,  Witten's formula (\ref{Wittenformula}) can be rewritten as following
 \beq\label{zhangformula}
\tilde{D}^2=\nabla^*\nabla+\frac{1}{4}(R+H^2-2e^0\cdot dH\cdot).
\eeq
Using (\ref{zhangformula})  we  studied lower bounds for the eigenvalues of the
Dirac-Witten operator under the modified energy condition or under
some other suitable conditions\cite{Chen2009}. We obtained, in the limiting case, the spacelike hypersurface is either maximal and Einstein manifold with positive
scalar curvature or Ricci-flat manifold with nonzero constant mean curvature.
We also refer to \cite{Maerten} for the case of bounded domains with smooth boundary and to \cite{CHZ2011,Chen2012}
for the case of higher codimensions.

In this paper,  we consider the locally reducible spacelike submanifold $(M^m,g)$ of  pseudo-Riemannian manifold $(\widetilde{M}^{m,n}, \tilde{g}). $ That is,
$TM=T_1\oplus\cdots\oplus T_k $ is orthogonal sum, where $T_i$ are parallel distributions of dimension $m_a,a=1,\cdots,k,$ and $m_1\geq m_2\geq\cdots\geq m_k.$
By introducing  appropriate modified connections, we give optimal
lower bounds' estimates for the eigenvalue of the Dirac-Witten
operator on spacelike submanifold $M^m\hookrightarrow\widetilde{M}^{m,n}.$  These estimates are given in terms of intrinsic and extrinsic expressions and
the limiting cases are examined. As a corollary, our estimate contains  some  results in \cite{A1,A,Chen2009,Chen2012} as special cases.
 The results can be considered the semi-Riemannian version of  our recent work in the eigenvalue estimates of the submanifold Dirac operator \cite{Chen2020}.

\mysection{\textbf{Preliminaries} } \ls
\subsection{Locally reducible Riemannian  manifolds} \ls

Let $M^m$ be a closed Riemannian spin manifold with nonnegative scalar curvature. Suppose
$TM=T_1\oplus\cdots\oplus T_k $ is orthogonal sum, where $T_i$ are parallel distributions of dimension $n_a,a=1,\cdots,k,$ and $n_1\geq n_2\geq\cdots\geq n_k.$ One important
consequence of the parallelism of $T_a$ is that $R(X,Y)=0$ whenever $X\in T_a, Y\in T_b$ with
$a\neq b.$
Then one can define
a locally decomposable
Riemannian structure $\beta$ as follows
\beq
\beta\mid_{T_1}=\textrm{Id},\ \  \beta\mid_{T_1^{\perp}}=-\textrm{Id}.
\eeq
Suppose
$$\{e_1,\cdots,e_{m_1},e_{m_1+1}, \cdots,e_{m_1+m_2 },\cdots,e_{m_1+m_2+\cdots+m_{k-1}+1 },\cdots, e_m\}$$
is an adapted local orthonormal frame, i.e., such that
\beQ
\{e_{m_1+m_2+\cdots+m_{a-1}+1},\cdots,e_{m_1+m_2+\cdots+m_a}\}
\eeQ
spans the subbundle $T_a.$
 Let $I_a=\{m_{a-1}+1,\cdots,m_a\},$
then
\beq
D_a&=&\sum_{i\in I_a} e^i\cdot\nabla_i
\eeq
is the ``partial" Dirac operator of subbundle $T_a,$ which  is formally self-adjoint operator.
Then
\beq
D&=&D_1+D_2+\cdots+D_k\\
D_\beta&=&D_1-(D_2+\cdots+D_k).
\eeq
Note for $a\neq b,$ one has
\beq
D_aD_b+D_aD_b=0,
\eeq
that is,
\beq\label{D2sum}
D_\beta^2=D^2=\sum_{a=1}^kD^2_a.
\eeq
In addition, we also have
\beQ
DD_\beta+D_\beta D=2(D_1^2-D_2^2-\cdots-D_k^2).
\eeQ

Let $R_a$ be  the ``scalar curvature" of $T_a,$ i.e.,
$$R_a=\sum_{s,t\in I_a}\langle R(e_s,e_t)e_t, e_s\rangle.$$
Hence the scalar curvature $R$ of $M^m$ is $R=\sum_{a=1}^kR_a.$
\subsection{The Dirac-Witten operator} \ls
\quad

Let $\widetilde{M}^{m,n}$ be an $(m+n)$-dimensional pseudo-Riemannian manifold whose metric
$\tilde{g}$ has signature
$(\underbrace{1,\cdots,1}_m,\underbrace{-1,\cdots,-1}_n).$ Let $M^m$ be an
$m$-dimensional spacelike submanifold with its induced Riemannian  metric $g$.
Let $\tilde{\nabla}$ and ${\nabla}$ be the Levi-Civita  connections of $\widetilde{M}^{m,n}$ and $M^m$ respectively.
Choose an orthonormal frame $\{e_\alpha\}$ and its dual
basis $\{e^\alpha\}$  with $e_i$ tangent  and $e_A$ normal to $M^m$
( throughout this paper, we agree on the following ranges of
indices:
$$1\leq\alpha,\beta,\gamma,\cdots\leq m+n;\ \ \ \ \ 1\leq i, j, k,\cdots\leq m$$
$$m+1\leq A, B, C,\cdots\leq m+n,$$
and the Einstein summation notation is also used).

 Suppose that
$M^m$ is a spin spacelike submanifold whose normal bundle  in $\widetilde{M}^{m,n}$ is also  spin.
Let $\textup{Spin}_0(m,n)$ be the connected component of  the pseudo-orthogonal spin group,
which is a  double cover of the connected component $\textup{SO}_0(m,n)$ of the pseudo-orthogonal group.
Denote by $K$ the
maximal compact  subgroup of  $\textup{Spin}_0(m,n)$ which covers
$\textup{SO}(m)\times \textup{SO}(n) \subset\textup{SO}_0(m,n).$

Let $\textit{\$}$ be the
(local) spinor {bundle} of $\widetilde{M}^{m,n}$.
  Since the following diagram is commutative
\[
\xymatrix{
&K \ar[rr]^{i}   \ar[d]_{\textup{Ad}}&  & \textup{Spin}_0(m,n)\ar[d]_{\textup{Ad}}\\
& \textup{SO}(m)\times \textup{SO}(n) \ar[rr]^{i}&& \textup{SO}_0(m,n) }
\]
the induced spinor bundle $\Bbb S\triangleq\textit{\$}|_M$ is globally defined over $M^m$.

Denote also by $\tilde{\nabla}$ and ${\nabla}$ the spin connections on $\Bbb S$. It is well-known \cite{Baum} that
 there exists a Hermitian inner product
$(\cdot,\cdot )$ on $\textit{\$}$  which is compatible with the
spin (local) connection $\tilde{\nabla}$. Moreover, for any vector field
 $\tilde{X}\in \Gamma(\widetilde{M}^{m,n})$
 and spinor fields $\phi,\psi\in\Gamma(\textit{\$})$, we have
\beq\label{nondefinite}
(\tilde{X}\cdot\phi,  \psi)=(-1)^{n+1}(\phi,  \tilde{X}\cdot\psi),
\eeq
where $``\cdot"$ denotes the Clifford multiplication with respect to the metric $\tilde{g}$.
Note that
this inner product is not positive definite.

The above diagram also implies that, over $M^m$,  there
exists on spinor bundle $\Bbb S$ a positive definite
Hermitian inner product defined by
$$\langle \cdot,\cdot \rangle:=(\omega\cdot, \cdot ),$$
where $\omega:=(\sqrt{-1})^{\frac{n(n-1)}{2}}e^{m+1}\cdots e^{m+n}$ (see \cite{Baum,CHZ2011}).
Obviously, $\omega^2=1$ and we have
$$\langle\omega\cdot\phi, \omega\cdot\psi  \rangle=\langle\phi, \psi  \rangle,$$
which implies
$$\langle\omega\cdot\phi, \psi  \rangle=\langle\phi, \omega\cdot\psi  \rangle.$$
\\
A simple computation yields
\begin{lemma}\label{lemma1}
Suppose that $M^m$ is spin submanifold whose  normal bundle is also spin. Over $M^m$, we have
\beq
\langle e^i\cdot\phi, \psi \rangle=-\langle\phi, e^i\cdot\psi  \rangle,\quad
\langle e^A\cdot\phi, \psi \rangle=\langle\phi, e^A\cdot\psi  \rangle.
\eeq
\end{lemma}
Fix a point $p\in M^m$ and an orthonormal basis
$\{e_\alpha\}$ of $T_p\widetilde{M}^{m,n}$ with $\{e_A\}$ normal and $\{e_i\}$ tangent to
$M^m$. Extend $\{e_i\}$ to a local orthonormal frame in a neighborhood
of $p$ in $M^m$ such that $(\nabla_ie_j)_p=0$. Then
\beq(\tilde{\nabla}_ie^A)_p=-h_{Aij}e^j+\nabla_i^\bot e^A,\quad\quad
(\tilde{\nabla}_ie^j)_p=-h_{Aij}e^A,
\eeq
where $h_{Aij}:=-\tilde{g}(\tilde{\nabla}_i e^A, e_j), 1\leq
i,j\leq m$, are the components of the second fundamental form and
$\nabla_i^\bot e^A=a_{iAB}e^B=-a_{iBA}e^B.$
Then one can deduce
 the following spinorial Gauss type formula:
\begin{equation}
\label{spinorgauss} \tilde{\nabla}_i=\nabla_i+\frac{1}{2}h_{Aij}e^A\cdot e^j.
\end{equation}

It is easy to check that,
\begin{lemma}\label{lemma2}
For any spinor field $\phi\in \Gamma(\Bbb S)$, we have
\beq
\nabla_i(e^A\cdot\phi)&=&(\nabla_i^\bot e^A)\cdot\phi+e^A\cdot\nabla_i\phi,\label{n}\\
\nabla_i(\omega\cdot\phi)&=&\omega\cdot\nabla_i\phi.\label{o}
\eeq
\end{lemma}

In the above orthonormal coframe $\{e^i\} $ of $M^m$, the
Dirac operator and the Dirac-Witten operator of $M^m$ are defined respectively by
$$D=e^i\cdot\nabla_i,\quad \tilde{D}=e^i\cdot\tilde{\nabla}_i.$$
From (\ref{spinorgauss}), they are related by
\begin{equation}
\label{2}
\tilde{D}=D+\frac{1}{2}\vec{H},
\end{equation}
where $\vec{H}=H_{A}e^A=\sum_i h_{Aii}e^A$ is the mean curvature covector field of $M^m$. Denote
$$| \vec{H}|:=\sqrt{\sum_AH_{A}^2}$$
throughout the paper.

It was proved \cite{CHZ2011}
that the connection $\nabla$ is compatible  with $\langle\cdot,\cdot\rangle.$
 and moreover, both $D$ and $\tilde{D}$ are formally
self-adjoint with respect to  $\int_M\langle\cdot,\cdot\rangle,$  if $M^m$ is closed.
In fact,
\beQ
\langle D\varphi,\psi\rangle=\langle e_i\cdot \nabla_i\varphi,\psi\rangle
&=&-\langle \nabla_i\varphi,e_i\cdot \psi\rangle\\
&=&-e_i\langle\varphi,e_i\cdot \psi\rangle+\langle \varphi,\nabla_i (e_i\cdot \psi)\rangle\\
&=&-e_ig(X,e_i)+\langle \varphi, D\psi\rangle\\
&=&-\textup{div}(X)+\langle \varphi, D\psi\rangle.
\eeQ
By Lemma \ref{lemma1} and (\ref{2}), it follows
\beQ
\langle \tilde{D}\varphi,\psi\rangle
&=&-\textup{div}(X)+\langle \varphi, \tilde{D}\psi\rangle.\\
\eeQ

\begin{proposition}
For any spinor field $\phi\in \Gamma(\Bbb S)$, we have
\beq\label{Chenformula}
\tilde{D}^2=\nabla^*\nabla+\frac{R+|\vec{H}|^2}{4}-\frac{1}{8}e^i\cdot e^j\cdot \textup{Id}\otimes\mathcal{R}_{ij}^\bot
 +\frac{1}{2}D^\bot\vec{H},
 \eeq
 where $R$ is the scalar curvature of $(M^n,g)$,
 $D^\bot\vec{H}:=e^i\cdot \nabla^\bot_i\vec{H}$  and $\mathcal{R}_{ij}^\bot$ stands for the
spinor normal curvature tensor.
\end{proposition}
\emph{{Proof}}.
It is straightforward to check that
\beq
\tilde{D}^2\phi
&=&\left(D+\frac{1}{2}\vec{H}\right)^2\phi\nonumber\\
&=&D^2\phi+\frac{1}{2}D(\vec{H}\cdot\phi)+\frac{1}{2}\vec{H}\cdot D\phi+
\frac{|\vec{H}|}{4}^2\phi\nonumber.
\eeq
But Lemma 2 implies for any spinor field $\phi$,
\beq \label{*}
(D^\bot\vec{H})\cdot\phi
=D(\vec{H}\cdot\phi)+\vec{H}\cdot D\phi.
\eeq
At the same time, the square $D^2$ is given by the following Schr\"{o}dinger-Lichnerowicz  formula
(see\cite{LM}, p. 164)
\beq
{D}^2=\nabla^*\nabla+\frac{R}{4}-\frac{1}{8}e^i\cdot e^j\cdot \textup{Id}\otimes\mathcal{R}_{ij}^\bot.
\eeq
\qed

 On the complement set  ${M_\phi}:=\{x\in M^m: \phi(x)\neq0\}$ of zeroes  of a spinor field
$\phi\in\Gamma(\Bbb S)$, define the functions
\beq
\check{R}_\phi&:=& R+|\vec{H}|^2+2\langle D^\bot\vec{H}\cdot\phi,\phi/|\phi|^2\rangle,\label{CR}\\
R^{\bot}_\phi&:=&-\frac{1}{2}\langle (e^i\cdot e^j\cdot \textup{Id}\otimes\mathcal{R}_{ij}^\bot)\phi,\phi/|\phi|^2\rangle.\label{RB}
\eeq
Hence, integrating (\ref{Chenformula}) on a closed manifold $M^m$, we obtain the following integrated version of
Schr\"{o}dinger-Lichnerowicz type formula
\beq\label{integratedChenformula}
\int_M |\tilde{D}\phi|^2=\int_M |\nabla \phi|^2+\frac{1}{4}(\check{R}_\phi+R^{\bot}_\phi)|\phi|^2
\eeq

\mysection{\textbf{Main results}} \ls

\begin{theorem}\label{theorem1}
Let $M^m\subset \widetilde{M}^{m,n}$ be a closed spacelike spin submanifold of dimension $m\geq 2,$ $TM^m=T_1\oplus\cdots\oplus T_k,$ where $T_a$ are parallel and pairwise orthogonal distributions of dimension $m_a, a=1,\cdots,
k, $ and
$m_1\geq m_2\geq\cdots\geq m_k, m_1>1.$ Let  $\lambda_H$ be any nonzero
eigenvalue of the Dirac-Witten operator $\tilde{D},$  and assume
$m_1(\check{R}_\psi+R^{\bot}_\psi)>-|\vec{H}|^2$
on $M_\psi,$
 then
\beq\label{main result H1}
\lambda_H^2\geq\frac{1}{4(m_1-1)^2}\inf_{M_\psi}
\left(\sqrt{m_1(m_1-1)(\check{R}_\psi+R^{\bot}_\psi)+m_1|\vec{H}|^2}-|\vec{H}|\right)^2.
\eeq
If $\lambda_H^2$ achieves its minimum and the normal bundle is flat, then both $R$ and $|\vec{H}|$
must be constant.
\end{theorem}

\proof
First, one can define
a locally decomposable
Riemannian structure $\beta$ as follows
\beQ
\beta\mid_{T_1}=\textrm{Id},\ \  \beta\mid_{T_1^{\perp}}=-\textrm{Id}.
\eeQ
 and  we also define the following modified connection, for any spinor field $\phi\in\Gamma(\Bbb S)$
\beq
T_i\phi
&=&\nabla_i\phi
+\left(\frac{p}{2}\vec{H}\cdot+q\lambda_H\right)\frac{1}{2}(\beta(e_i)+e_i)\cdot\phi+
\frac{1}{2}\nabla_{(\beta-\textrm{Id})(e_i)}\phi,
\eeq
where $p, q$ are smooth real-valued functions that are specified later.
If $\tilde{D}\psi=\lambda_H\psi$ for a nontrivial $\psi,$ a direct  computation gives
\beq
|T\psi|^2\label{|T|}
&=&|\nabla \psi|^2+(m_1pq-p-q)\mathfrak{R}e\langle \vec{H}\cdot D\psi,\psi\rangle+(m_1q^2-2q)|\tilde{D}\psi|^2\nonumber\\
&{}&+\frac{|\vec{H}|^2}{4}(m_1p^2+2q-2m_1pq)|\psi|^2-\sum_{i=m_1+1}^m|\nabla_i\psi|^2\nonumber\\
&&-\frac{1}{2}\mathfrak{R}e\langle (p\vec{H}\cdot-2q\lambda_H)\psi,(D-D_\beta)\psi\rangle.
\eeq
By the Weitzenb\"{o}ck type formula (\ref{integratedChenformula}), integrating (\ref{|T|}) over $M^m$ yields
\beq\label{integral term}
&{}&\int_M|T\psi|^2+\sum_{s=m_1+1}^m|\nabla_i\psi|^2+\frac{1}{2}\mathfrak{R}e\langle (p\vec{H}\cdot-2q\lambda_H)\psi,(D-D_\beta)\psi\rangle\nonumber\\
&=&\int_M\left(m_1pq-p-q\right)\mathfrak{R}e\langle \vec{H}\cdot D\psi,\psi\rangle+(m_1q^2-2q+1)|\tilde{D}\psi|^2\nonumber\\
&{}&+\frac{|\vec{H}|^2}{4}(m_1p^2+2q-2m_1pq)|\psi|^2-\frac{\check{R}_\psi+R^{\bot}_\psi}{4}|\psi|^2.
\eeq
Let $m_1pq-p-q=0,$ where
\beq\label{q}
(m_1q-1)^2=\frac{(m_1-1)|\vec{H}|}{\sqrt{m_1(m_1-1)(\check{R}_\psi+R^{\bot}_\psi)+m_1|\vec{H}|^2}-|\vec{H}|}.
\eeq
Note if $q=1/m_1,$ one has $ \vec{H}=0,$ and in this case, $\tilde{D}=D.$ Therefore, it would be well that we also assume  $q\neq1/m_1,$ in the following discussion.

From $(\ref{q}), $ one has the following
\beq\label{IRHS}
&&\textrm{r.h.s. of (\ref{integral term})}\nonumber\\
&=&\int_M\left(m_1q^2-2q+1\right)\nonumber\\
&{}&\times\left[\lambda_H^2-\frac{1}{4}\left(\frac{\check{R}_\psi+R^{\bot}_\psi}{m_1q^2-2q+1}
-\frac{m_1p^2+2q-2m_1pq}{m_1q^2-2q+1}|\vec{H}|^2\right)\right]|\psi|^2\nonumber\\
&=&\int_M\left(m_1q^2-2q+1\right)\nonumber\\
&{}&\times\left[\lambda_H^2-\frac{1}{4(m_1-1)^2}
\left(\sqrt{m_1(m_1-1)(\check{R}_\psi+R^{\bot}_\psi)+m_1|\vec{H}|^2}-|\vec{H}|\right)^2\right]|\psi|^2.
\quad\quad\eeq

Now we turn to deal with the l.h.s. of (\ref{integral term}).
First, observe that for any spinor field $\phi\in \Gamma(\mathbb{S})$
\beq\label{T}
\sum_{i=1}^me_i\cdot T_i\phi=D\phi+m_1\left(\frac{1}{2}p\vec{H}\cdot-q\lambda_H\right) \phi+\frac{1}{2}(D_\beta\phi-D\phi).
\eeq
Therefore, with the help of $D^2=D^2_\beta$
\beq\label{key}
&{}&\mathfrak{R}e\int_M\langle \sum_{i=1}^me_i\cdot T_i\phi, \frac{1}{2}(D\phi-D_\beta\phi)\rangle\nonumber\\
&{=}&\frac{ m_1}{4}\int_M\mathfrak{R}e\left\langle \left(p\vec{H}\cdot-2q\lambda_H\right) \phi,D\phi-D_\beta\phi\right\rangle.
\eeq
 Hence, using the Cauchy-Schwarz inequality and the equality (\ref{D2sum}) we can deduce that
\beq\label{ILHS}
&{}&\textrm{l.h.s. of}\ (\ref{integral term})\nonumber\\
&=&\int_M|T\psi|^2+\sum_{s=m_1+1}^m|\nabla_i\psi|^2+\frac{1}{2}\mathfrak{R}e\langle (p\vec{H}\cdot-2q\lambda_H)\psi,(D-D_\beta)\psi\rangle\nonumber\\
&=&\int_M|T\psi|^2+\sum_{s=m_1+1}^m|\nabla_s\psi|^2-\frac{2}{m_1}\mathfrak{R}e\left\langle \sum_{i=1}^{m_1}e_i\cdot T_i\psi, \frac{1}{2}(D_\beta\psi-D\psi)\right\rangle\nonumber\\
&\geq&\int_M|T\psi|^2+\sum_{s=m_1+1}^m|\nabla_s\psi|^2-\frac{1}{\sqrt{m_1}}\left(\varepsilon^{-1}| T\psi|^2+\frac{\varepsilon}{4}|D_\beta\psi-D\psi|^2\right)\nonumber\\
&\geq&\int_M\left(\frac{1}{m_2}-\frac{1}{m_1}\right)|D_2\psi|^2+\cdots+\left(\frac{1}{m_k}
-\frac{1}{m_1}\right)|D_k\psi|^2,
\eeq
here we take $\varepsilon=\frac{1}{\sqrt{m_1}}.$

Hence, (\ref{ILHS}),  combined with   (\ref{IRHS}),  yields that
\beQ
0&{\leq}&\int_M
\left(\frac{1}{m_2}-\frac{1}{m_1}\right)|D_2\psi|^2
+\cdots+\left(\frac{1}{m_k}
-\frac{1}{m_1}\right)|D_k\psi|^2\nonumber\\
&\leq&\int_M(m_1q^2-2q+1)\nonumber\\
&&\times\left[\lambda_H^2-\frac{1}{4(m_1-1)^2}
\left(\sqrt{m_1(m_1-1)(\check{R}_\psi+R^{\bot}_\psi)+m_1|\vec{H}|^2}-|\vec{H}|\right)^2\right]|\psi|^2.\quad
\eeQ
Therefore, the first part of the theorem follows.

 If $\lambda_H^2$ achieves its minimum, then
\beq
\sum_{i\in I_1}e_i\cdot T_i\psi=\frac{1}{2}(D_\beta\psi-D\psi),
\eeq
 which, together with (\ref{T}), implies that,
\beq\label{D1}
D\psi+m_1\left(\frac{1}{2}p\vec{H}\cdot-q\lambda_H\right)\psi=0.
\eeq
From the definition of $\tilde{D}$ and  the relationship  of smooth functions $p$ and $q,$
one obtains
\beq
\vec{H}\cdot\psi
&\stackrel{(\ref{D1})}{=}&2(1-m_1q)^2\lambda_H\psi\label{H1}\\
&\stackrel{(\ref{q})}{=}&\textup{sign}(\lambda_H)|\vec{H}|\psi.\label{H2}
\eeq
Moreover,
\beq\label{D2}
D\psi=f\psi,
\eeq
 where $f:=\lambda_H (2 m_1 q-m_1^2q^2).$

Now suppose $m_1=m_2=\cdots=m_l>m_{l+1}\geq\cdots \geq m_k, $  where $l\geq 1,$ then
\beq
D_{l+1}\psi=D_{l+2}\psi=\cdots= D_k\psi=0
\eeq
and we also have, for $i,j\in I_1,$ the two spinor fields are proportional, i.e.,
\beq
e_i\cdot T_i\psi=e_j\cdot T_j\psi.
\eeq
And also for any $a\in \{1,\cdots,l\},$ and any $s,t\in I_a,$
\beq
e_s\cdot \nabla_s\psi=e_t\cdot \nabla_t\psi.
\eeq
and $\alpha\in (I_1\cup I_2\cup\cdots\cup I_l)^c,$
\beq
\nabla_\alpha\psi=0.
\eeq
Therefore, for $i\in I_1,$
\beq
\nabla_i\psi
&=&T_i\psi-\left(\frac{1}{2}p\vec{H}\cdot+q\lambda_H\right) e_i\cdot\psi\nonumber\\
&\stackrel{(\ref{D1})}{=}&-\frac{1}{2m_1}e_i\cdot (D_\beta\psi-D\psi)- \frac{1}{m_1} e_i\cdot D\psi\nonumber\\
&=&-\frac{1}{2m_1}e_i\cdot D\psi-\frac{1}{2m_1}e_i\cdot D_\beta\psi\nonumber\\
&\stackrel{(\ref{D2})}{=}&-\frac{1}{2m_1}fe_i\cdot \psi-\frac{1}{2m_1}e_i\cdot D_\beta\psi
\eeq
and while, for any $b\in \{2,\cdots,l\},s\in I_b,$
\beq
\nabla_s\psi
&=&-\frac{1}{m_1}e_s\cdot D_b\psi\nonumber\\
&=& \frac{1}{2m_1}e_s\cdot (D_\beta\psi-D\psi)+\frac{1}{m_1}\sum_{a=2,\neq b}^l e_s\cdot D_a \psi\nonumber\\
&\stackrel{(\ref{D2})}{=}&- \frac{1}{2m_1}fe_s\cdot\psi+\frac{1}{2m_1}e_s\cdot D_\beta\psi
+\frac{1}{m_1}\sum_{a=2,\neq b}^l e_s\cdot D_a \psi.
\eeq
Furthermore,  for  any $b\in \{2,\cdots,l\}$ and any $s,t\in I_b,$ one can check the following
\beq
\nabla_t\nabla_s\psi
&=&-\frac{f_t}{2m_1}e_s\cdot\psi+\frac{f^2}{4m_1^2}e_s\cdot e_t\cdot\psi
-\frac{f}{4m_1^2}e_s\cdot e_t\cdot D_\beta\psi-\frac{ f}{2m_1^2}e_s\cdot e_t\cdot \sum_{a=2,\neq b}^l  D_a \psi\nonumber\\
&{}&+\frac{1}{2m_1}e_s\cdot\nabla_t(D_\beta\psi)+
\frac{1}{m_1}\sum_{a=2,\neq b}^l e_s\cdot \nabla_t(D_a \psi).
\eeq
Since the normal bundle is flat, it follows that
\beq
-\frac{1}{2}Ric(e_t)\cdot\psi
&=&\sum_{s\in I_b}e_s\cdot\mathcal{R}_{e_t,e_s}\psi
\nonumber\\
&=&\sum_{s\in I_b}e_s\cdot[\nabla_t,\nabla_s]\psi\nonumber\\
&=&\frac{1}{2m_1}(m_1f_t+D_bf\cdot e_t)\cdot \psi-\frac{m_1-1}{2m_1^2}f^2e_t\cdot\psi
+\frac{m_1-1}{2m_1^2}fe_t\cdot D_\beta\psi\nonumber\\
&{}&+\frac{m_1-1}{m_1^2}fe_t\cdot\sum_{a=2,\neq b}^l  D_a\psi+\frac{1}{2m_1}[(-m_1+2)\nabla_tD_\beta\psi+e_t\cdot D_b(D_\beta\psi)]
\nonumber\\
&{}&+\frac{1}{m_1}\left[(-m_1+2)\nabla_t\left(\sum_{a=2,\neq b}^l  D_a\psi\right)+e_t\cdot D_b\left(\sum_{a=2,\neq b}^l  D_a\psi\right)\right].
 \eeq
Eventually, we obtain for any $b\in \{2,\cdots,l\}$
\beq
\frac{1}{2}R_b\psi
&=&-\frac{1}{2}\sum_{t\in I_b} e_t\cdot Ric(e_t)\cdot\psi\nonumber\\
&=&\frac{m_1-1}{m_1}D_b(f)\cdot \psi+\frac{m_1-1}{2m_1}f^2\psi\nonumber\\
&{}&-\frac{m_1-1}{2m_1}f D_\beta\psi-\frac{m_1-1}{m_1}f \left(\sum_{a=2,\neq b}^l D_a\psi\right)\nonumber\\
&{}&-\frac{m_1-1}{m_1}D_b(D_\beta\psi)-\frac{2(m_1-1)}{m_1}D_b\left(\sum_{a=2,\neq b}^l D_a\psi\right).\ \ \ \ \ \ \ \
\eeq
And while, the same argument gives rise to
\beq
\frac{1}{2}R_1\psi
&=&-\frac{1}{2}\sum_{i\in I_1} e_i\cdot Ric(e_i)\cdot\psi\nonumber\\
&=&\frac{m_1-1}{m_1}D_1(f)\cdot \psi+\frac{m_1-1}{2m_1}f^2\psi\nonumber\\
&{}&+\frac{m_1-1}{2m_1}f D_\beta\psi+\frac{m_1-1}{m_1}D_1(D_\beta\psi).
\eeq
Taking sum, one obtains
\beq
\frac{1}{2}R\psi
&=&\frac{1}{2}\sum_{a=1}^l R_a \psi\nonumber\\
&=&\frac{m_1-1}{m_1}\left(\sum_{a=1}^lD_af\right)\cdot \psi+l\cdot\frac{m_1-1}{2m_1}f^2\psi\nonumber\\
&{}&-(l-2)\cdot\frac{m_1-1}{2m_1}f D_\beta\psi+
\frac{m_1-1}{m_1}\left(D_1-\sum_{b=2}^lD_b\right)(D_\beta\psi)\nonumber\\
&{}&-(l-2)\cdot\frac{m_1-1}{m_1}f\sum_{b=2}^lD_b\psi.
\eeq
Note, for $\alpha\in (I_1\cup I_2\cup\cdots\cup I_l)^c,$ we have
$\nabla_\alpha\psi=0.$
Hence
\beq
\sum_{b=2}^lD_b\psi=\frac{1}{2}(D-D_\beta)\psi
\eeq
and
\beq
\left(D_1-\sum_{b=2}^lD_b\right)(D_\beta\psi)
&=&D_\beta^2\psi\nonumber\\
&=&D^2\psi\nonumber\\
&=&D(f\psi)\nonumber\\
&=&\nabla f\cdot\psi+f^2\psi.
\eeq
Hence
\beq
\frac{1}{2}R\psi&=&\frac{2(m_1-1)}{m_1}\left(\sum_{a=1}^lD_af\right) \cdot \psi+\frac{2(m_1-1)}{m_1}f^2 \psi\nonumber\\
&{}&\quad+\frac{m_1-1}{m_1}\left(\nabla f-\sum_{a=1}^lD_af\right)\cdot\psi,
\eeq
which implies that $\nabla f=0$ and $f^2=\frac{m_1}{4(m_1-1)}R.$
The whole proof of Theorem 1 is complete.
\qed

\begin{corollary}\label{corollary}
Let $M^m\subset \widetilde{M}^{m,1}$ be a closed spacelike spin hypersurface of dimension $m\geq 2,$ $TM^m=T_1\oplus\cdots\oplus T_k,$ where $T_a$ are parallel and pairwise orthogonal distributions of dimension $m_a, a=1,\cdots,
k, $ and
$m_1\geq m_2\geq\cdots\geq m_k, m_1>1.$ Let  $\lambda_H$ be any nonzero
eigenvalue of the Dirac-Witten operator $\tilde{D},$  and assume
$m_1\check{R}_\psi>-H^2$
on $M_\psi,$
 then
\beq\label{main result H2}
\lambda_H^2\geq\frac{1}{4(m_1-1)^2}\inf_{M_\psi}
\left(\sqrt{m_1(m_1-1)\check{R}_\psi+m_1H^2}-|H|\right)^2.
\eeq
If $\lambda_H^2$ achieves its minimum,  there are only two possibilities in the limiting case:\\
(i):  $H\equiv0$ and
 \beq
 \lambda^2_H\equiv \frac{m_1}{4(m_1-1)} R.
 \eeq
The universal covering space  of $M^m$
is isometric to a product $M_1\times \cdots\times M_k,$ where $dim M_s=m_s, M_1$ has a real Killing
spinor and  $M_s$ has a parallel spinor if $m_s<m_1,$  and $M_s$ has a parallel spinor or a real
 Killing spinor if $m_s=m_1.$\\
(ii): $M^m$ is Ricci-flat, and \beq\lambda^2\equiv \frac{H^2}{4}.\eeq
\end{corollary}
\proof
Since we already know $f,H$ and $ R $ are all constant in the limiting case, it follows that
\beq
0=(D^\bot \vec{H})\cdot\psi
&=&D(\vec{H}\cdot\psi)+\vec{H}\cdot D\psi\nonumber\\
&=&2(m_1q-1)^2\lambda_HD\psi+f\vec{H}\cdot \psi\nonumber\\
&=&4(m_1q-1)^2\lambda_Hf\psi\nonumber\\
&=&4(\lambda_H-f)f\psi.
\eeq
which yields $f=0$
(recall the previous assumption that  $q\neq1/m_1$).
 In this case, (\ref{q}) implies that the scalar curvature $R=0$ on the whole manifold and the parallelism of $T_a$ implies that  $R_a, a=1,\ldots,k,$ are  constant.

Note, for any $k\in I_a, a\in \{1,\cdots,l\},$
\beq
\nabla_k\psi+\frac{1}{m_a}e_k\cdot D_a\psi=0.
\eeq
and $\alpha\in (I_1\cup I_2\cup\cdots\cup I_l)^c,$
\beq
\nabla_\alpha\psi=0.
\eeq
Hence, following the arguments in \cite{A1}, one knows that in the limiting case,  the universal covering space of $M^m$
is isometric to a product $M_1\times \cdots\times M_k,$ where $\dim M_s=m_s, s=1,\cdots,k, $
each $M_s$ has a parallel spinor.
The whole proof of corollary \ref{corollary} is complete.


\begin{thebibliography} {99}

\bibitem{A1}B. Alexandrov, The first eigenvalue of the Dirac operator on locally reducible
           Riemannian manifolds, J. Geom. Phys. 57 (2007), no. 2, 467--472.
\bibitem{A} B. Alexandrov, G. Grantcharov, S. Ivanov, An estimate for the first eigenvalue
         of the Dirac operator on compact Riemannian spin manifold admitting a parallel
          one-form, J. Geom. Phys. 28 (1998), no. 3-4, 263--270.
\bibitem{Baum}  H. Baum,  A remark on the spectrum of  the  Dirac operator on pseudo-Riemannian spin manifold (1996) (preprint).

\bibitem{Bar93}  C. B\"ar, Real Killing spinors and holonomy, Comm. Math. Phys. { 154} (1993), 509--521.
\bibitem{CHZ2011} D. Chen, O. Hijazi, X. Zhang,   The Dirac-Witten operator on
         pseudo-Riemannian manifolds, Math. Z. 271 (2012), 357--372.

\bibitem{Chen2020}Y. Chen,	Lower bounds for the eigenvalue estimates of the submanifold Dirac operator, arXiv:2010.13016.

\bibitem{Chen2009} Y. Chen,   Lower bounds for eigenvalues of the Dirac-Witten operator,
                 Sci. China. Ser. A-Math.  { 52} (2009), 2459--2468.

\bibitem{Chen2012} Y. Chen, Xu Xu,  Some remarks on the Dirac--Witten operator on pseudo-Riemannian manifolds, J. Geom. Phys.  62 (2012), 1999--2008.



\bibitem{Bar93}  C. B\"ar, Real Killing spinors and holonomy, Comm. Math. Phys. { 154} (1993), 509--521.

\bibitem{Friedrich}T. Friedrich, Dirac operators in Riemannian geometry,
        Graduate Studies in Mathematics 25, American Mathematical Society, 2000.

\bibitem{HZ03} O. Hijazi, X. Zhang,   The Dirac-Witten operator on
         spacelike hypersurfaces,  Comm. Anal. Geom. { 11} (2003), 737--750.

\bibitem{LM} H. B. Lawson,  M. L. Michelsohn, Spin Geometry,
         Princeton  Math Series. {38} Princeton University Press, 1989.
\bibitem{Maerten}D. Maerten, Optimal eigenvalue estimate for the Dirac--Witten operator
on bounded domains with smooth boundary. Lett. Math. Phys. 86
(2008), 1--18.

\bibitem{PT} T.  Park, C.  Taubes,  On Witten's proof of the
         positive energy theorem, {Comm. Math. Phys.} { 84} (1982),  223--238.

\bibitem{Witten} E. Witten,   A new proof of the positive energy theorem,
          Comm. Math. Phys.  { 80} (1981), 381--402.


\bibitem{Yano} K. Yano, M. Kon,  Structures on manifolds. Singapore: World Sci. 1984.

\bibitem{ZZ} L. Zhang,  X. Zhang, { Remarks on  positive energy theorem,}
             Comm. Math. Phys. { 208} (2000), 663--669.

\bibitem{Zhang1999} X. Zhang,  Positive mass theorem for modified energy
         condition, proceedings of the workshop on Morse theory, minimax
          theory and  their applications to  nonlinear differential
          equations held at Morningside Center of Mathematics, (1999), 275--282.





\end{thebibliography}
\end{document}